\documentclass{article}
\usepackage{amsmath,amssymb,enumerate,graphicx}
\newtheorem{thm}{Theorem}
\newtheorem{mydef}{Definition}
\newtheorem{cor}{Corollary}

\newtheorem{prop}{Proposition}
\newtheorem{assum}{Assumption}

\newcommand{\R}{\mathbb R}
\newcommand{\A}{{\cal{A}}}
\newcommand{\N}{{\cal{N}}}

\DeclareMathOperator{\viab}{viab}
\newcommand\be{\begin{equation}}
\newcommand\nd{\end{equation}}
\usepackage{caption}
\usepackage{subcaption}
\usepackage[section]{placeins}
\usepackage{float}
\usepackage{multirow}

\title{Collective management of environmental commons with multiple usages: a guaranteed viability  approach.}
\author{{\sc I. Alvarez${}^{1,}{}^{2}$, L. Zaleski${}^{3}$, J-P. Briot${}^{3}$, M. de A. Irving${}^{4}$ }\\[2mm]
	${}^{1}$ UCA, INRAE, LISC \\
9 avenue Blaise Pascal CS 20085, F-63178  Aubi\`{e}re, France\\
e-mail: {\tt Isabelle.Alvarez@inrae.fr}\\	
	${}^{2}$ CNRS, ISC-PIF\\
113 rue Nationale, F-75013, Paris, France\\[2mm]
	${}^{3}$ SU, CNRS, LIP6\\
	4 place Jussieu, F-75005, Paris, France\\
	e-mail: {\tt laetitia.zaleski@gmail.com}, {\tt Jean-Pierre.Briot@lip6.fr}\\[2mm]
		${}^{4}$ EICOS, IP, UFRJ\\
		Avenida Pasteur, 250 Rio de Janeiro 22290-240, BR\\
	e-mail: {\tt mirving@mandic.com.br}
}

\date{}

\begin{document}
\maketitle

\section*{Abstract}
In this paper we address the collective management of environmental commons with multiple usages in the framework of the mathematical viability theory. 
We consider that the stakeholders can derive from the study of their own socioeconomic problem the variables describing their different usages of the commons and its evolution, and a representation  of the desirable states for the commons. We then consider the guaranteed viability kernel, subset of the set of desirable states where it is possible to maintain the state of the commons even when its evolution is represented by several conflicting models. 
This approach is illustrated on a problem of lake eutrophication.

\textit{Key words} viability theory; guaranteed viability; collaborative decision; environmental commons; lake eutophication

\section{Introduction}
Sustainable use of natural resources, environmental conservation, social inclusion and welfare, economic activity and development imply generally conflicting management objectives. In \textit{the tragedy of the commons}, \cite{hardin1968} highlights the exhaustion of open-access resources by numerous users with similar view, but the same analysis can be done with different types of users whose activity is based on the resource.
A lot of work on sustainability of natural resources is still focused on allocation problem, where stakeholders are considered as competitors for the share of quotas, for example for the regulation of fisheries or water sharing (see references for instance in \cite{parrachino_cooperative_2006}, \cite{zaccour_survey_2018}).

In order to take into account the different types of stakeholders' interests, many efforts are done in the economic approach to assess the value of environmental and social services (see for example a framework in \cite{deGroot2005}). When points of view are considered to be incommensurable, multi-criteria or viability theory approaches propose interesting alternatives. Even when stakeholders are considered as competitors for one common resource, these approaches allow to take into account more indicators than the level of renewable resource and the profit directly based on it. For example, in a quantitative work on fishing regulation \cite{DOWLING2020109243}, 21 score functions are designed for the regulation of fishing, depending on fish biomass and parameters computed each year depending on the control scenario. Weighted stakeholders' preferences over the score functions are then optimized each year for different levels of the control variable.
When stakeholders express different points of view, the viability theory (VT) approach allows to combine the different constraints placed on the system, without direct connection to the underlying profit of the related activity. For instance, in a hydro-power dam management problem \cite{Alais_DeLara_multi-usage_2017}, the main concern is maximizing the profit of the electricity provider with water control under uncertainty on water inflow and electricity price. Recreational and agricultural activities impose an additional seasonal constraint on the water level without further profit analysis. In \cite{wei_sustainability_2013}, the multi-objective concern of a tourist city is studied through the linked evolution of the number of tourists, tourism infrastructure and environment quality. The different stakes are represented by constraints on the level of these variables. The VT algorithm identifies the area where it is possible to maintain the evolution of the three variables between these bounds.
 Viability approach has shown its potential in many other domains as stated in the review from \cite{zaccour_survey_2018}.
In all these works, the model of the evolution of common resources or land uses, together with the impact of controls on the system (such as the total allowable catch in fishery regulation), is supposed to be consensual. It is generally taken from the literature or from previous work and parameters are calibrated from data and time series.  
In works cited in \cite{zaccour_survey_2018}, the model representing the system at stake is always considered as consensual. In theoretical works, models can certainly be generic functions of variables, controls and uncertainties (it is the case in \cite{23_de_lara_multi-criteria_2009},\cite{53_krivan_differential_1991}, \cite{55_krivan_non-stochastic_1998}, \cite{45_martinet_risk_2016}). But in their applications and in the other works cited in \cite{zaccour_survey_2018}, when uncertainty is explicitly taken into account, it is in fact related to data and measurements used to assess parameters in the model, not to the definition of the model itself. As stated in \cite{45_martinet_risk_2016}, uncertainty affects parameters such as growth rate, recruitment or mortality in dynamic population models, unknown or unpredictable events such as climate fluctuations, or externalities such as price, as in games against nature. Models are supposed to be consensual with their explicit hypotheses (which are generally discussed).
Actually in \cite{44_1_little_elfsimmodel_2007} three models are considered for larval dispersal and modelers can parameterize the system to run simulations with their own choice of model. It is motivated by the possibility of studying different species, so for the viability study only one model is parameterized. 

However, the ComMod approach \cite{ComMod_2014} has shown that modeling the evolution of the system at stake is difficult and hardly consensual, since scientific or technical viewpoint can be considered by stakeholders as a viewpoint among others. 
The ComMod approach addresses this problem with serious games supported by simulation models \cite{o_barreteau_role-playing_nodate} where stakeholders can test their hypotheses about the system evolution and the impact of actions. The process goes on, with additional research if necessary, until a consensus on the model is reached. 
 
To take into account this discrepancy at the model level we consider here that stakeholders have their own model of the evolution of the system with the impact of controls. We consider that stakeholders are able to define constraints on the key variables of their usages of the commons, such as the number of tourists, the quality of water (for example measured in term of concentration of pollutant or bacteria), the quality of the environment (for example measured with biodiversity indicators related to the population of local species), etc. These constraints are generally seen as thresholds.
We consider that the objective of the group of stakeholders is to define a set of states where the system can be maintained with appropriate controls. We use viability theory, in particular the concept of guaranteed viability set \cite{AubinGuaranti}, which is  defined to take into account uncertainties (such as move by nature, see for instance \cite{BATES2018244}).   

The paper is organized as followed: we first describe the problem and our hypotheses, together with a reminder of viability theory. We describe individual and group viewpoint as viability problems, and show why in this context of several models it is more difficult to seek technically sound agreement. We illustrate this approach with a problem of lake eutrophication. In Section 3, we first present the perturbation embedding function, which allows us to consider each model as a perturbation of a central model. This formulation enables the definition of a guaranteed viability problem. 
In Section 4 we present and discuss the application to the management of the lake. We summarize the results and perspectives in the concluding section.

\section{Definition of the problem}
\label{sec:defprob}
Let us consider an entity $\A$ (for instance a preserved area) in evolution, which is submitted to the management decision of a group of $N$ members. 
The group has to first define a project for $\A$ as a set $K$ of desirable states within which the state of $\A$ should remain with an appropriate set of actions depending on the state of $\A$. The group has then to find a solution to this project. 

Each member $i\in {\cal{N}}:=\{1,\dots ,N\}$ has a personal view and project for entity $\A$. Each member $i$ considers that entity $\A$ is described by a vector of state variables $x_i\in  \mathbb{R}^{n_i}$, and that its 
evolution is governed by a controlled dynamical system. Member $i$ considers that the state of $\A$ should remain in a set of desirable states, $K_i \subset \mathbb{R}^{n_i}$. Member $i$ deems that the controls which should be used depend on the state of $\A$. These admissible controls are defined by a set-valued map $U_i:\,\mathbb{R}^{n_i}\leadsto \mathbb{R}^{p_i}$, where $U_i(x_i)$ is the set of admissible controls that member $i$ finds appropriate at state $x_i$. 
\begin{mydef}
	\label{def:project_i} $(K_i,U_i)$ defines member $i$'s project for the management of $\A$, where $K_i$ is the set of desirable sets and $U_i$ the admissible control map.
\end{mydef}
From the viability theory viewpoint, a solution to member $i$'s project consists in finding the states of $K_i$ from which there always exists a control function selected in $U_i$ so that the state of $\A$ remains in $K_i$.


The objective of the group to define its project, compatible with every member's project, then to find a solution to it, compatible with everyone's solution.


We use as an illustration a problem of lake eutrophication as stated in \cite{carpenter_management_1999}. Agricultural practice and other human activities can lead to lake pollution with phosphates. Phosphorus dynamics in the lake can lead to eutrophication, which negatively impacts the biodiversity of ecosystem, and causes serious annoyance to residents and tourism activities. We consider that a committee is formed to study and manage the problem. It is composed of farmers and local elected authorities.

\subsection{What members share about their project}
\label{sec:CommonHypotheses}

\begin{assum}
\label{A1}
Members of the decision group share the knowledge about which state variables they consider for $\A$, so $\forall i \in {\cal{N}}, n_i=n$. We note $x \in \mathbb{R}^{n}$ the vector of state variables of $\A$ that is shared by each member of the group.
\end{assum}
In particular, when measurements of the state of $\A$ are possible, all members agree on the validity of the measure, so if the measure of the state of $\A$ at date $T$ is $x(T)$, then $\forall i \in {\cal{N}}, x_i(T)=x(T)$.

In the lake and nearby farms problem, following \cite{carpenter_management_1999} and \cite{Martin2004}, we consider that all members agree that the key variables to the problem are the Phosphorus input (noted $L$) and the Phosphorus concentration in the lake (noted $P$). 
\begin{assum}
\label{A2}
Members of the decision group share the knowledge about the control variables they consider for $\A$. We note $U: \mathbb{R}^{n}\leadsto \mathbb{R}^{p}$ the set-valued map of admissible control which associates the group's set of admissible controls with the state of $\A$: $\forall i \in {\cal{N}}, U_i=U$.
\end{assum}
Assumption \ref{A2} supposes that all group member's have reached a consensus regarding which control variables are to be considered. This can be achieved by restriction to the intersection of members' control set so that $\forall x \in \mathbb{R}^{n}, U(x)=\bigcap_{i\in {\cal{N}}} U_i(x)$. When the latter is empty, a negotiation should take place to build a non empty $U(x)$. We suppose here that such a negotiation has taken place.

In the lake and nearby farms problem, following \cite{Martin2004}, we consider that the committee agrees to the possibility of controlling the rate of variation of the Phosphorus input and to maintain this rate between boundaries, so $U=\left[ u_{min},u_{max}\right]$. This can be done by farmers controlling their fertilizer input, by the greater or lesser use of wetlands or by the use of water treatment plants (\cite{gajardo_modeling_2017}).

With Assumption \ref{A1} and \ref{A2}, it is possible to define the group project for $\A$.
\begin{mydef}
\label{def:project}
Let $K=\bigcap_{i\in {\cal{N}}} K_i$ with $K \neq \emptyset$, and $U: \, \mathbb{R}^{n}\leadsto \mathbb{R}^{p}$ a set-valued control map defined on $K$. $(K,U)$ is a group project if all members agree that the state of $\A$ should remain in the set of desirable states $K$, 
using controls from the admissible map $U$.
\end{mydef}
When $K \subset \mathbb{R}^{n} = \emptyset $, a negotiation should take place to build a non empty $K$. We suppose here that such a negotiation has taken place.

In the lake and nearby farms problem, everybody wants to keep the lake in an oligotrophic state, which supposes to set a concentration limit of Phophorus $P_{max}$ in the lake (for example established from previous observations). Everybody also wants to maintain or develop the agricultural activity, which supposes to allow a minimum amount of Phophorus input $L_{min} $ in the lake. So everybody agree to maintain the state of the lake described by $(L,P)$ in a set of desirable states $ K=\left[ L_{min},+\infty\right)  \times \left[ 0,P_{max}\right]$.

When Assumption \ref{A1} and Assumption \ref{A2} are verified, all group members can describe the dynamics of the state of $\A$ as they seen it as a (possibly discrete) controlled dynamical system with the shared variables. 
In the case of the lake and neighboring farms, the dynamics for member $i$ is supposed to be defined according to \cite{carpenter_management_1999} and \cite{Martin2004}, by the equations of system \ref{eq:lake}  with the constraints on $L_{min}$ and $P_{max}$.
\begin{equation}
	S(b_i,r_i,q_i,m_i)\left\{
	\begin{array}{lll}
		\frac{dL}{dt}&=&u\in U=\left[ u_{min},u_{max}\right] \\
		\frac{dP}{dt}&=&-b_i P(t) + L(t) +r_i\frac{P(t)^{q_i}}{m_i^{q_i} + P(t)^{q_i}}\\
	\end{array}\right.
	\label{eq:lake}
\end{equation}
$u_{min}$ and $u_{max}$ are the maximum effort farmers and local authorities are ready to allow or to take to increase or decrease the phosphorus input. 
$b_i$ is the rate of loss (due to sedimentation and outflow),  $r_i$, $q_i$ and $m_i$ are parameters of the sygmoid-like (s-shaped) dynamics of Phophorus recycling in the lake, which are generally set by calibration from observations: $r_i$ is the maximum rate of recycled Phophorus, $m_i$ is the concentration of Phosphorus at which the recycling rate is half its maximum and $q_i$ is a parameter of the steepness of the dynamics (see \cite{carpenter_management_1999} for more details).

In the general case, we note $Sc(f,U)$ the continuous dynamical system defined by:
	\begin{equation}
Sc(f,U)\left\{
\begin{array}{lcl}
x'(t)&=&f(x(t),u(t))\\
u(t)&\in &U(x(t))\subset \mathbb{R}^{p}\; .
\end{array}
\right.
\label{evolve}
\end{equation}
where $f$ is a function from $ \mathbb{R}^{n}\times \mathbb{R}^{p}$ to $ \mathbb{R}^{n}$ and $U$ a set-valued map from $\mathbb{R}^{n}$ to $ \mathbb{R}^{p}$. Similarly, we note $Sd(f,U)$ the discrete dynamical system defined by:
	\begin{equation}
Sd(f,U)\left\{
\begin{array}{lcl}
x^{k+1}&=&f(x^k,u^k)\\
u^k&\in &U(x^k)\subset \mathbb{R}^{p}\; .
\end{array}
\right.
\label{evolveD}
\end{equation}
We note $Sc_i=Sc(f_i,U)$ the dynamical system that described the evolution of $\A$ for member $i$ in a continuous case. We note $Sd_i=Sd(f_i,U)$ the discrete dynamical system that described the evolution of $\A$ for member $i$ in a discrete case. The function $f_i: \mathbb{R}^{n}\times \mathbb{R}^{p}\rightarrow \mathbb{R}^{n}$ associates the variations of $\A$ state variables for member $i$ with the current values of the state and control variables.   

We don't assume that the different stakeholders in the group share their dynamics. We consider here that they don't necessarily agree on dynamics, and that they are not compelled to make their belief public. But we assume that they agree to share this information with a trusted third party.  

\subsection{A reminder of the viability theory}
\label{sec:Viability}
Referring to \cite{Aubin91}, we define viable evolutions and the viability kernel. 
\begin{mydef} An evolution of the system $Sc(f,U)$ (\ref{evolve}) (resp. $Sd(f,U)$ (\ref{evolveD})) is viable in $K$ 
	if and only if its trajectory remains in $K$. In the continuous case: $\forall t\in \mathbb{R}^+ \  x(t)\in K.$ In the discrete case: $\forall k\in \mathbb{N} \  x^k\in K.$
	\label{def:viable_evolv}
\end{mydef}
\begin{mydef} A set $L$ is viable for the system $Sc(f,U)$ (\ref{evolve}) (resp. $Sd(f,U)$ (\ref{evolveD})) if for all $x \in L$ there is an evolution of the system $Sc(f,U)$ (\ref{evolve}) (resp. $Sd(f,U)$ (\ref{evolveD})) starting at $x$ and viable in $L$.
		\label{def:viable_set}
\end{mydef}

\begin{mydef} \label{def:kernel} The viability kernel associated to system $Sc(f,U)$ (\ref{evolve}) (resp. system $Sd(f,U)$ (\ref{evolveD})) under constraint $K$ is the set of all states in $K$ from which there is an evolution of $Sc(f,U)$ (resp. $Sd(f,U)$) starting at $x$ and viable in $K$.
\end{mydef}

Under some general conditions listed in appendix \ref{annexe1}, the viability kernel is a close set. In the interior of the viability kernel, all control are viable, so viable controls on the boundary show how it is possible to maintain the system in the constraint set. This information can be used to define control strategies.

\begin{mydef}
A control map with images restricted to viable controls only is called a viable regulation map.
\label{def:viablereg}
\end{mydef} 

\begin{prop}
\cite{Aubin91}. If $L$ is a viable set for the system $Sc(f,U)$ (\ref{evolve}) (resp. $Sd(f,U)$ (\ref{evolveD})), let $\tilde{U}$ be a viable regulation map, then $L$ is a viable set for the system $Sc(f,\tilde{U})$ (\ref{evolve}) (resp. $Sd(f,\tilde{U})$ (\ref{evolveD}). Moreover, for all $x \in L$, any evolution starting from $x$ and governed by $Sc(f,\tilde{U})$ (resp. $Sd(f,\tilde{U})$ in the discrete case) is viable in $L$. $L$ is called an invariant set for dynamics $Sc(f_i,\tilde{U})$ (resp. $Sd(f_i,\tilde{U})$.
	\label{prop:invariantSet}
\end{prop}
From any state in the viability kernel, it is always possible to find a control function that allows the state of the system to stay in the viability kernel indefinitely. Conversely, from any initial state outside the viability kernel, there is no way to prevent the exit in finite time of an evolution governed by system (\ref{evolve}) (resp. (\ref{evolveD}) in the discrete case). 

\begin{figure}[t!]
	\centering
	\includegraphics[width=8cm]{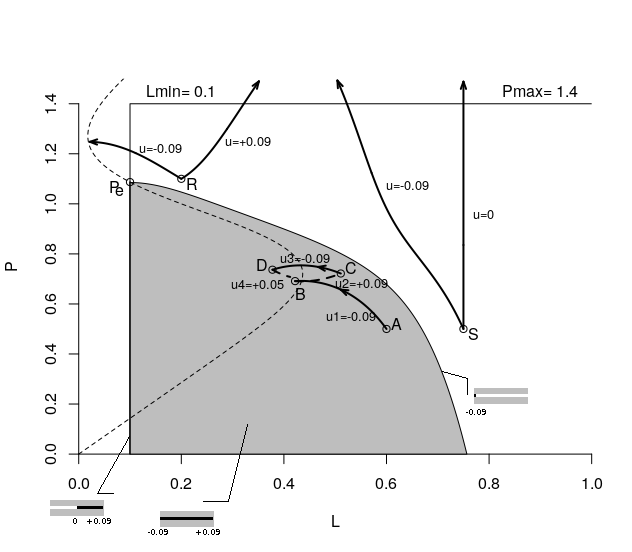}
	\caption{\small{Viability kernel (in gray) of the lake and neighboring farms problem, with $L_{min}=0.1$, $P_{max}=1.4$ ($L$ and $P$ in $\mu$gL$^{-1}$), $U=\left[ -0.9,0.9\right]$, dynamics parameters value $q=8$, $m=r=1$, $b=0.7$. Constraint set boundary is in black plain lines ($L=L_{min}$, $P=P_{max}$). The curve of equilibria is dashed. Viable controls are shown as a black line in cartouches. A viable trajectory starting from A is shown ($u=-0.09$ from $A$ to $B$, then a cycle with $u=+0.09$ from $B$ to $C$, $u=-0.09$ from $C$ to $D$, $u=+0.05$ from $D$ to $B$). From S where the lake is still in an oligotrophic state, even with maximum effort from the farmers ($u=u_{min}$) the concentration of Phosphorus becomes too high. From state R also outside the viability kernel, all trajectories leave the constraint set, leading the lake to eutrophic state or farmers' activity to an unsustainable state.}}
	\label{fig:trajlac}
\end{figure}

In the case of the lake and its neighboring farms, it is shown in \cite{Martin2004} that the viability kernel associated to system (\ref{eq:lake}) submitted to the constraint $(L,P)\in K=\left[ L_{min},+\infty\right)  \times \left[ 0,P_{max}\right]$ is not empty when $P_{max}$ is greater than the smallest $P$-value of the equilibria associated with $L_{min}$ (an equilibrium $P$-value is defined by $\frac{dP}{dt}=0$).
For example in Figure \ref{fig:trajlac}, the state $(L_{min},P_e)$ is an equilibrium with $P_e \leq P_{max}$, so the viability kernel is not empty. When the curve of Equilibria intersects the half-line $(L\geq L_{min},P=P_{max})$ at $(L_e,P_{max})$, the boundary of the viability kernel is delimited by the segment line $(L=L_{min},P\leq P_{max})$, the segment line $(L_{min} \leq L\leq L_{e},P=P_{max})$ and the integral curve of the dynamics with control $u=u_{min}$ arriving in $(L_e,P_{max})$.  When the curve of Equilibria doesn't intersect the half-line $(L\geq L_{min},P=P_{max})$, as in Figure \ref{fig:trajlac}, we note $P_e$ the P-value of the highest equilibrum on the segment line $(L=L_{min},P\leq P_{max})$. In that case, the boundary of the viability kernel is delimited by the segment line $(L=L_{min},P\leq P_{e})$ and the integral curve of the dynamics with control $u=u_{min}$ passing through $(L_{min},P_{e})$. 

Figure \ref{fig:trajlac} shows the viability kernel for the lake and neighboring farm problem in this latter case, for a given set of parameters for system~(\ref{eq:lake}) and constraint set $K$. From any state in this viability kernel, it is possible to find a trajectory that stays in the viability kernel indefinitely. Figure \ref{fig:trajlac} presents an example of viable trajectory from a state in the viability kernel. It also shows examples of states outside the viability kernel; even the most severe control cannot prevent trajectories from leaving the constraint set. Either the lake will shift to an eutrophic state, or the economic activity will be jeopardized. 

From a state outside the viability kernel, every evolution governed by system~(\ref{eq:lake}) with this particular set of parameters, choice of constraint set and control interval will exit the constraint set. In general, dealing with state outside the viability kernel implies to study the resilience (as in \cite{Martin2004}) or to redefine the problem. This can be done by relaxing the constraints on the desirable set (when it is possible),  by allowing more efficient control which are not presently part of the admissible controls, or by modifying the dynamics. This latter option is generally more difficult to implement, since it implies to modify the lake itself (see \cite{LakeChinois} for example of such actions).

\subsection{Viewpoints as viability problems}
\label{sec:viewpoints}
\paragraph{Objectives.}
\label{sec:Objectives}

\begin{figure}[t!]
	\centering
	\includegraphics[height=6cm]{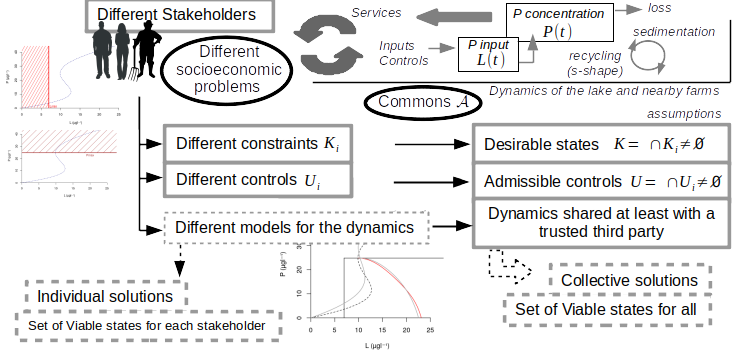}
	\caption{\small{Diagram of the finding of management solutions for the lake and nearby farms (system $\A$) in the framework of viability theory with different stakeholders. Gray arrows denote relationship in the dynamics model. Large gray arrows represents interaction with stakeholders models and dynamics (which are not explicit). Black arrows represents the viability analysis process. Dotted lines and arrows show the main focus of the article. }}
	\label{fig:diagram}
\end{figure}

We assume that Assumptions \ref{A1} and \ref{A2} are verified, and that the group has defined $(K,U)$ as its project for $\A$ according to Definition \ref{def:project}. $K \neq \emptyset \subset \mathbb{R}^{n}$ is the set of desirable states for $\A$ and $U: \mathbb{R}^{n}\leadsto \mathbb{R}^{p}$ the set-valued map of admissible controls. With this definition of the group project, each member can work on a solution according to the dynamics he assumes for $\A$. 
Finally, the objective of the group is to design a solution from all member's solutions. Figure \ref{fig:diagram} summarizes the implications of considering different usages and stakeholders for the lake and nearby farms problem. Although the intuition is to work from the set of individual solutions, in this section we show that this approach is difficult to implement. 

\subsubsection{Individual Viewpoint}
\label{sec:individual}
We consider here that each member $i$ is able to describe the evolution of the state of $\A$ with a controlled dynamical system. It is either a continuous system $Sc_i=Sc(f_i,U)$ (equation \ref{evolve}), or a discrete system $Sd_i=Sd(f_i,U)$ (equation \ref{evolveD}).
We also consider in the following that the conditions for Proposition \ref{prop:kernelclosed_c} (resp. \ref{prop:kernelclosed_d}) are fulfilled: the viability kernel associated to member $i$'s system and constraints is closed.
%

Let $L_i\subset K$ be a non-empty viable set for the continuous (resp. discrete) system $Sc(f_i,U)$ (\ref{evolve}) (resp. $Sd(f_i,U)$ (\ref{evolveD})) submitted to constraints $K$. Then from all states in $L_i$ there is at least one viable evolution governed by  $Sc(f_i,U)$ (resp. $Sd(f_i,U)$) that stays in $L_i$. From member $i$ viewpoint, $L_i$ is a solution state set to the management of $\A$.
\begin{mydef}
	$L_i\subset K$ is a solution state set for Member $i$ for project $(K,U)$ if $L_i$ is a non-empty viable set for Member $i$'s dynamics.
	\label{def:solutioni}
\end{mydef}

We note $\viab_i(K)$ the viability kernel associated to member $i$'s project with dynamical system $Sc(f_i,U)$ (\ref{evolve}) (resp. $Sd(f_i,U)$ (\ref{evolveD})) submitted to the viability constraint $K$. 
In the case of the lake and its neighboring farms, Figure \ref{fig:trajlac} shows the viability kernel for the dynamics (\ref{eq:lake}) submitted to constraint set $K=\left[ L_{min},+\infty\right)  \times \left[ 0,P_{max}\right]$ for the particular values of the dynamics parameters (noted as farmers representative in Figure \ref{fig:celac}).

\subsubsection{Group Viewpoint}
In the following we suppose that the group project for  $\A$ is $(K,U)$ and that all group members can propose their own individual solution to the management of  $\A$: $$\forall i \in {\cal{N}}, \viab_i(K) \neq \emptyset$$

We note $H=\bigcap_{i\in {\cal{N}}} \viab_i(K)$. If $H=\emptyset$ a negotiation should obviously take place between stakeholders, since there is no way to operate $\A$ and satisfy the group members. When the intersection is not empty, it seems a good candidate.  
The intersection of viability kernels has already been proposed as a solution to insure the viability of two fishing fleets operating on the same resource. In \cite{38_sanogo_viability_2012}, the intersection of the viability kernel of both fleets is viable for each fleet if they change their effort at the same time when necessary, which suppose a high level of cooperation. 
But unfortunately, it is not always the case.

\begin{prop}
The intersection of all members' viability kernels is not necessary viable for all members. 
\label{prop:admiNOTconsensus}
\end{prop}

\begin{figure}[b!]
	\centering
\includegraphics[height=6cm]{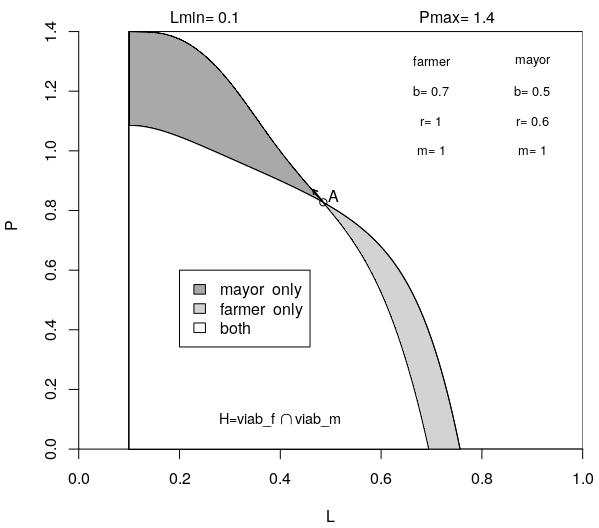}
\caption{\small{Viability kernels of two stakeholders in the lake and neighboring farms problem, with $L_{min}=0.1$ and $P_{max}=1.4$ ($L$ and $P$ in $\mu$gL$^{-1}$), $U=\left[ -0.9,0.9\right]$, shared parameters value $q=8,m=1$ . In white, the intersection $H$ of the viability kernels. In dark (resp. light) gray, the complementary area of the mayor (resp. farmer) viability kernel. State A is not viable in the intersection for the mayor. The arrow shows the trajectory of state A according to the mayor: it leaves the white area.}}
\label{fig:celac}
\end{figure}

Proof. The problem of the lake and neighboring farms gives a counter-example. We consider here two stakeholders, say a mayor and a farmer's union representative (respectively noted with $m$ and $f$ indices). Both stakeholders interpret the observations in different ways, so they adopt different values for the parameters of dynamics (\ref{eq:lake}). Their respective viability kernel ($\viab_m$ and $\viab_f$) associated to the constraint set $K=\left[ L_{min},+\infty\right)  \times \left[ 0,P_{max}\right]$  are shown on Figure \ref{fig:celac}. For the particular parameters chosen, the intersection $H$ is not empty. For $x=(L,P) \in viab_i$, viable controls are defined by $\tilde{U_i}$ with $\tilde{U_i}(L_{min},P)=[0,u_{max}]$, $\tilde{U_i}(L,P)=\{u_{min}\}$  when $(L,P)$ is on the boundary of $\viab_i$ with $L\neq L_{min}$, and otherwise $\tilde{U_i}(x)=U$.  Nevertheless, state $A$ in the intersection is not viable for the mayor. State $A$ is on the boundary of the mayor's viability kernel, so $\tilde{U}_m(A) =\{u_{min}\}$ and the only viable control for the mayor is $u=-0.09$. But the trajectory starting at $A$ and governed by $(S_m)$ with $u=-0.09$ stays on the boundary of $\viab_m$ so it leaves $H$. $\blacksquare$

\begin{mydef}
	\label{def:interU}
	Let $L\subset K$ and let $(U_i)_{i\in {\cal{N}}}$ be control maps defined on $L$. Let $U$ be defined for all $x\in L$ by $U(x)= \bigcap_{i\in {\cal{N}}} U_i(x)$. 
$U$ is called the intersection of $(U_i)_{i\in {\cal{N}}}$ on $L$.  
\end{mydef} 	
We note  $\tilde{U}$ the regulation map defined on the intersection $H$ of all viability kernels of the group members by the intersection of all the corresponding viable regulation map: $\tilde{U}(x)=\bigcap_{i\in {\cal{N}}} \tilde{U}_i(x)$. 
Obviously, if there is a state $z\in H$ such that $\tilde{U}(z)=\emptyset$, it means that members cannot agree on a way to control the evolution of $\A$ at this particular state. Unfortunately, even if all members agree on controls on $H$, it is not sufficient to reach a consensus.
\begin{cor}
	\label{cor:domNOTsuff}
Let $\tilde{U}$ be the intersection on $H=\bigcap_{i\in {\cal{N}}} \viab_i(K)$ of the viable regulation map $U_i$ on each $ \viab_i(K)$ of each member $i\in \N$. $Dom(\tilde{U}) =H$ is not a sufficient condition for $H$ being a viable set for all members.
\end{cor}
Proof. In the previous example, we can derive that $\tilde{U}_m(L,P)=\tilde{U}_f(L,P)$ for all $(L,P)$ in the intersection except on the set $H_b$ of the boundary of $H$ wherever $L\neq L_{min}$. On the part of the boundary of $H$ which is the boundary of $viab_m$ only, $\tilde{U_m}(L,P)=\{u_{min}\}$, while $\tilde{U_f}(L,P)=U$ (and conversely on the boundary of  $viab_f$ only). So for $(L,P) \in H_b$, $\tilde{U}(L,P) =\{u_{min}\}$ so $Dom(\tilde{U}) =H$. Nevertheless state $A$ is not viable in $H$ for the mayor. $\blacksquare$


From Proposition \ref{prop:admiNOTconsensus}, and Corollary \ref{cor:domNOTsuff} we propose the following definition for a technically sound consensus solution to the management of $\A$.
\begin{mydef}
	Let $(K,U)$ be the projet of the management group for $\A$. A set of state $H\subset K$ is a consensus solution if $H$ is a viable set for the each member and if the domain of the intersection $\tilde{U}$  of the viable regulation maps of each member on $H$, $(U_i)_{i\in {\cal{N}}}$, is such that $Dom(\tilde{U})=H$.
	\label{consensus}
\end{mydef}

Actually, to be a consensus solution, a subset $H$ of the constraint set has to be viable for all members and each viable state has to share at least one viable control for all members. In that case it is possible for the group member to reach an agreement on the control, regardless of trajectories. For example, in the discrete case, from any state $x_0$ of $H$, all group members share at least one viable control value that allows the state of $\A$ to stay in $H$. Since the dynamics they consider are different, there is generally no consensus on state $x_1$. But as long as the group members still share a viable control value they still can agree on it. When it is no longer the case, for example at step $n$, the true value of the state of $\A$ can be measured to continue this process from $x_n$ as new starting point. In the continuous case, when $H$ is a close set, such situations arise only on the boundary of $H$.

\begin{figure}[t!]
	\centering
	\includegraphics[width=8cm]{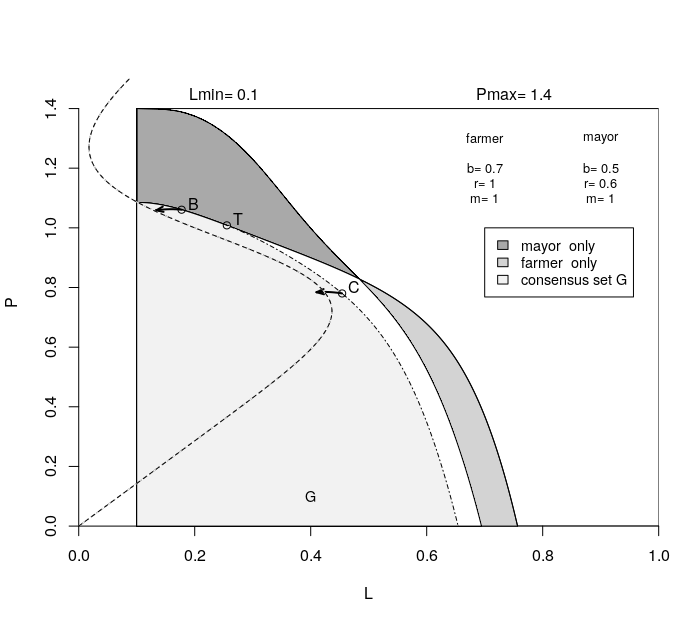}
	\caption{\small{A consensus state set for the two stakeholders in the lake and neighboring farms problem with parameters from Figure \ref{fig:celac}. In very light gray, the consensus state set G, delimited in dot dashed line by the trajectory following the mayor's dynamics that is tangent to the boundary of the viability kernel of the farmer at state T (in plain line). From each state $(L,P)$ of this trajectory before the tangent state T (with $L>L_T$), the evolution governed by the farmer's system starting at these states with $u=u_{min}$  leaves the boundary to evolve inside G, as it is shown for state C. Respectively, from each state of the boundary of the viability kernel of the farmer with $L_{min}<L<L_T$, the evolution governed by the mayor's system starting at these states with $u=u_{min}$ leaves the boundary to evolve inside $G$, as it is shown for state B. In dashed line, the line of equilibrium for the farmers representative's dynamics. }}
	\label{fig:celacgaranti}
\end{figure}

Figure \ref{fig:celacgaranti} shows a consensus state set for the lake and neighboring farms system (\ref{eq:lake}) with parameters of Figure \ref{fig:celac}. The associated regulation map $\tilde{U}$ is such that $\tilde{U}(L_{min},P)=[0,u_{max}]$, $\tilde{U}(L,P)=\{u_{min}\}$  when $(L,P)$ is on the boundary with $L\neq L_{min}$, and otherwise $\tilde{U}(L,P)=[u_{min},u_{max}]$. For every state $(L_0,P_0)$ of the consensus state set $G$, there is an evolution governed by system (\ref{eq:lake}) for each stakeholder, starting at $(L_0,P_0)$, with the same $u(0) \in \tilde{U}(L_0,P_0)$ that stays in $G$. In the interior of the intersection of the viability kernels this property is also verified since for every states in the interior all controls are viable in the viability kernel of each stakeholder. For states on the boundary with $L=L_{min}$, for both the mayor and the farmers' representative several evolutions are viable, in particular with $u=0$. The consensus state space is delimited by the boundary of the farmers'representative and the trajectory governed by the mayor's dynamics with the minimum control value that stays in the viability kernel of the farmers'representative with largest $L$ when $P=0$.

In the case of the lake and neighboring farms problem, with only two group members, it is possible to define a consensus state set because of the properties of the dynamics, for which the line of equilibria is known, and the viability kernels and the trajectories corresponding to minimum control value $u_{min}$ can be easily defined and computed \cite{Martin2004}).

For more general cases it is necessary to propose a method that can be applied without such knowledge. We present such a method in the next section.

\section{Consensus with guaranteed Viability}

\subsection{Embedding function for the dynamics}
\label{sec:PEF}

Since the group members have their own definition for the dynamics, all members can see others' definitions as perturbations of their own. We show here that is possible to define the dynamics of $\A$ embedding all members' definitions seen as perturbations. The dynamics of $\A$ depends on the state of $\A$, $x(t)$, on the control chosen in $U(x(t))$ and on perturbations occurring from a set $V(x(t))\subset \R^q$ that depends on the state of $\A$. In the continuous case we have:

\begin{equation}
Svc(f,U,V)\left\{
\begin{array}{lcl}
x'(t)&=&f(x(t),u(t),v(t))\\
u(t)&\in &U(x(t))\\
v(t)&\in &V(x(t)) \; 
\end{array}
\right.
\label{evolveV}
\end{equation}
In the discrete case:
	\begin{equation}
Svd(f,U,V)\left\{
\begin{array}{lcl}
x^{k+1}&=&f(x^k,u^k,v^k)\\
u^k&\in &U(x^k)\subset \mathbb{R}^{p}\\
v^k&\in &V(x^k)\subset \mathbb{R}^{q}\; ,
\end{array}
\right.
\label{evolveDV}
\end{equation}
where $f$ associates the new state of $\A$ with its present state, a control chosen in $U(x(t))$ and a perturbation in $V(x(t))$. $Scv(f,U,V)$ and $Sdv(f,U,V)$ are called dynamical controlled tychastic systems \cite{AubinGuaranti}.
\begin{mydef}
We say that System $Svc(f,U,V)$ (\ref{evolveV}) (resp. $Svd(f,U,V)$ (\ref{evolveDV})) embeds System $Sc(f_i,U)$(\ref{evolve}) (resp. $Sd(f_i,U)$ (\ref{evolveD})) for $i \in \N$, and call the corresponding pair $(f,V)$ an embedding solution if and only if:
\begin{equation}
\forall x \in K, \forall u \in U(x), \forall i \in \N, \exists v_{i,u,x} \in V(x), f_i(x,u)=f(x,u,v_{i,u,x})
\label{condV}
\end{equation}
\end{mydef}
We show in appendix \ref{annexe2} that under some general conditions a System (\ref{evolveV}) (resp. \ref{evolveDV} in the discrete case) can embed $Sc(f_i,U)$(\ref{evolve}) (resp. $Sd(f_i,U)$ (\ref{evolveD})) for all $i \in \N$.

For example, for the problem of the lake and neighboring farms, the dynamics for every group member are defined from $S_{lake}$ in equation (\ref{eq:lake}) by $f_i : \mathbb{R}^{2}\times \R\leadsto \mathbb{R}^{2}$, with:

\begin{equation}
f_i((x_1,x_2),u)=\left( 
\begin{array}{l}
u\\
-b_i x_2 + x_1 +r_i\frac{x_2^{q}}{m + x_2^{q}} \;
\end{array}
\right)
\label{ex1}
\end{equation}
where parameters $m$ and $q$ have consensus values among the group, whether parameters $b_i$ and $r_i$ have not. Then, by defining $V=[\min_{i\in \N}(b_i),\max_{i\in \N}(b_i)]\times[\min_{i\in \N}(r_i),\max_{i\in \N}(r_i)]$ and $f$ as:

\begin{equation}
f((x_1,x_2),u,v)=\left( 
\begin{array}{l}
u\\
v_1 x_2 + x_1 + v_2 \frac{x_2^{q}}{m + x_2^{q}} \;
\end{array}
\right)
\label{ex1G}
\end{equation}
with $v=(v_1,v_2) \in V$, equation (\ref{condV}) is verified, since in this simple case we have:
$$\forall x \in K, \forall u \in U(x), \forall i \in \N, v_{i,u,x} = (b_i,r_i).$$

In the following, we assume that the group has defined a map $f$ and a perturbation map $V$ such that System (\ref{evolve}) (resp. (\ref{evolveD}) in the discrete case) describes the dynamics of $\A$, embedding the viewpoint of all group members. 

The objective of the group is then to find a consensus solution for $\A$ which will guarantee the viability for each member with shared viable controls.

\subsection{Guaranteed Viability with embedding dynamics}
\label{sec:GV}
We recall here some definitions and properties of the mathematical theory of viability, from \cite{Aubin91} and \cite{Lavallee2020v2}, relatively to guaranteed viability.

\begin{mydef}
	A solution $x(.)$ of system (\ref{evolveV}) (resp.\ref{evolveDV}) is an evolution $(t \mapsto x(t))$ (resp. $(x^k)_{k\in \mathbb{N}}$) such that there is a control function $(t \mapsto u(t))$ (resp. $(u^k)_{k\in \mathbb{N}}$) and a perturbation function $(t \mapsto v(t))$ (resp. $(v^k)_{k\in \mathbb{N}}$) such that system (\ref{evolveV}) (resp.\ref{evolveDV}) is verified for almost all $t\geq \R^+$ (resp. for all $k\in \mathbb{N}$).
\end{mydef}

\begin{mydef}
	An evolution $x(.)$ (resp. $(x^k)$) solution of system (\ref{evolveV}) (resp.\ref{evolveDV}) is viable in $L$ if and only if its trajectory remains in $L$.
		\label{viableEvolv}
\end{mydef}

Following \cite{Aubin91}, \cite{doyen_lipschitz_2000} and \cite{Lavallee2020v2}, we recall the property of guaranteed viability. 

\begin{mydef}(From \cite{AubinGuaranti})
	A set $L$ verifies the property of guaranteed viability for $Svc(f,U,V)$ (\ref{evolveV}) (resp. $Svd (f,U,V)$ (\ref{evolveDV})) if there is a regulation map $\tilde{U}$ defined on $L$ with non empty subset of $U$ images, i.e.  $\forall x \in L, \tilde{U}(x)\neq \emptyset$ and  $\tilde{U}(x)\subset U(x)$ such that for all $x_0$ in $L$, all evolutions starting at $x_0$ and governed by $Svc(f,\tilde{U},V)$ (resp. $Svd(f,\tilde{U},V)$) are viable in $L$. 
	\label{guarantiSet}
\end{mydef}

\begin{mydef}
	The guaranteed viability kernel associated to a set $K$ is the largest set in $K$ with property of guaranteed viability (for $\lambda$-Lipschitz controls in the continuous case - see appendix \ref{annexe1} for definition).
\end{mydef}
We note $Guar_{Sdv(f,U,V)}$ the guaranteed viability kernel associated to a set $K$ for the discrete dynamics (\ref{evolveDV}). In the continuous case (\ref{evolveV}), we note it $Guar_{\lambda, Scv(f,U,V)}$.

We have seen that the intersection of each member's solution is not necessarily a solution for all members. We are going to show that the guaranteed viability kernel is a consensus solution (as in Definition \ref{consensus}).
 
Let $Scv(f,U,V)$ (resp. $Sdv(f,U,V)$ ) be an embedding solution for all group members, which fulfilled conditions of Proposition \ref{prop:viabguaranty}. Let $L\neq \emptyset$ be the guaranteed viability kernel for system $Scv(f,U,V)$ with $\lambda$ Lipschitz constant (resp. $Sdv(f,U,V)$) associated to constraint set $K$. Then we have the following property:
\begin{thm}
	The guaranteed viability kernel associated to $K$ for $Scv(f,U,V)$ (with $\lambda$ Lipschitz constant in the continuous case) (resp. $Sdv(f,U,V)$ in the discrete case ) is a consensus solution to the management of $\A$.
\label{th:viab_consensus}
\end{thm}
The demonstration can be found in appendix \ref{annexe1.2}. We note $\text{Guar}_K$ the guaranteed viability kernel and $\tilde{U}$ the associated viable regulation map. The basic idea is that in $\text{Guar}_K$, from member $i$'s perspective, an evolution governed by system $Sc'_i=Sc(f_i,\tilde{U})$ (resp. $Sd'_i=Sd(f_i,\tilde{U})$ in the discrete case) is also governed by  the embedding system $Scv(f,U,V)$ (resp. $Sdv(f,U,V)$), so it remains in $\text{Guar}_K$, therefore $\text{Guar}_K$ is a viable set for each member $i$. $\blacksquare$ 

Under some general conditions, the guaranteed viability kernel is a close set (see Proposition \ref{prop:viabguaranty} in appendix \ref{annexe1.2}).
If the guaranteed viability kernel is closed, it is possible to retrieve the value of viable controls on its boundary. This information could be used to anticipate or to design control strategies that keep an evolution away from the boundary. 

%

\section{Application to the problem of lake eutrophication}
\label{sec:appli}
\subsection{Lake Bourget case}
We consider Lake Bourget, which is the biggest lake located entirely within France. It is monitored by the inter-syndicate committee CISALP, which is in charge of the design, animation and management of contractual actions for depollution and restoration of Lake Bourget. The lake had experienced a long eutrophication period, since in 1974, incoming amount of $P$ in the lake was around 300 tons per year, in 1989 the in-lake concentration was above 150 $mg.m^{-3}$ (\cite{vincon-leite_modelisation_1990}) where OECD norms assess the in-lake concentration to a maximum of $P_{o}=10$ $mg.m^{-3}$ (equivalent of 36 tons) for the oligotrophic state (from \cite{Vollenweider}). Similarly a threshold  for mesotrophic state $P_{m}=35$ $mg.m^{-3}$ can be defined. With concentration above $P_m$ the lake is supposed to be in an eutrophic state. Since the lake is monitored and the data are available, it is possible to calibrate the equations of system \ref{eq:lake} for lake Bourget. For Phosphorus unit in $mg.m^{-3}$ states values and parameters $r$ and $m$ are divided by the volume of the lake in billions of $m^{3}$. The volume of the lake is $v=3.6 \; 10^9$ $m^{3}$ and supposed to be constant. Calibration coefficients from \cite{brias2018} are given in Table \ref{tab:lake}. 
\begin{table}
\caption {Parameters of the Lake Bourget model} \label {tab:lake} \centering
\begin{tabular}{lcccc}
	\hline
	Parameter & b & r & q & m \\
	\hline
	\multirow{4}{*}{Value} & \multicolumn{4}{c}{state unit in tons} \\
		  & 2.2676 & 367.04  & 2.222 & 96.85 \\
	\cline{2-5}	
	& \multicolumn{4}{c}{state unit in $mg.m^{-3}$ (or $\mu g.l^{-1}$)} \\
	  & 2.2676 & 101.96 & 2.222 & 26.90 \\
	\hline		
\end{tabular}
\end{table}

Lake Bourget offers multiple services apart from being a freshwater reserve. It is an area of major ecological interest for flora, fauna and the diversity of its biotopes. Several areas of the lake are classified as protected area. It supports a lot of tourism and recreational activities (water-sport, fishing, beaches, marina) and cultural activities linked in particular to historical heritage and literature. Several other services are being considered, such as the production of hydrothermal energy. 
Although the state of the lake has considerably improved, it is still considered as oligo-mesotrophic. Its dynamics can be unstable due to P loading and several blooms of cyanobacteria have been observed lately.

Agriculture is now the main source of $P$ loading since major prevention measures have been taken since 1980. In particular the effluents of water treatment plants are no longer rejected in the lake. 

The control of incoming $P$ is considered as essential because of the potential lagged impact of $P$ release from sediments (\cite{Jacquet2017}).

We consider a scenario where the CISALP in its mission of negotiation would approach agricultural unions, local representatives and managers of tourism activities, to form a committee in order to control $P$ loading in the lake to prevent eutrophication and its consequences. The effort on $P$ loading we consider is a limit on its variation as it is common in environmental actions. It can be implemented by changes in agricultural practice and by the use of wetlands or retention basins. Actually retention basins are being built to regulate the incoming of polluted water.

As in Section \ref{sec:defprob}, we consider that the committee members are aware of models regarding lake Bourget (such as \cite{brias2018}), but they can disagree on the value of parameters or even on the model formulation. We consider that they can agree on a set $K$ of desirable states and on a set of admissible control $U$, or at least they can consider several scenarios for the definition of theses sets.

\subsection{Scenarios and Results}
\label{sec:results}

We consider a scenario where members of the committee agree on the state variables and the possibility of controlling the rate of P loading ($L$). They consider different parameters sets for model (\ref{eq:lake}), and possibly a different formulation for the process of recycling from sediments. Some members consider that the recycling process can have more effect at low value of $P$ total than with model (\ref{eq:lake}). A different formula for the sigmoid-like function is used in that case as shown in Equation (\ref{eq:Bourget2}), with parameter $\lambda_i$ controlling the shape instead of $q_i$ in model (\ref{eq:lake}). For small value of $\lambda_i>0$ the recycling occurs also for low level of in-lake $P$, so the lower branch of the equilibrium curve is actually higher.
\begin{equation}
S'_i \left\{
\begin{array}{lll}
\frac{dL}{dt}&=&u\in U=\left[ u_{min},u_{max}\right] \\
\frac{dP}{dt}&=&-b_i P(t) + L(t) + r_i \frac{P(t)}{P(t) + m_i e^{(-\lambda_i(P(t)-m_i)) }}
\end{array} \right.
\label{eq:Bourget2}
\end{equation}
It is possible to embed both model types by considering an additional parameter $\alpha_i \in [0,1]$ which controls the predominance of one type over the other. The corresponding model is represented in Equation (\ref{eq:lakeBourget}).

\begin{equation}
B_i \left\{
\begin{array}{lll}
\frac{dL}{dt}=u\in U=\left[ u_{min},u_{max}\right] \\
\frac{dP}{dt}=-b_i P(t) + L(t) +(1-\alpha_i)  r_i\frac{P(t)^{q_i}}{m_i^{q_i} + P(t)^{q_i}} + \alpha_i  r_i \frac{P(t)}{P(t) + m_i e^{(-\lambda_i(P(t)-m_i)) }}
\end{array} \right.
\label{eq:lakeBourget}
\end{equation}

\begin{figure}[!b]
	\centering
	\begin{subfigure}[b]{0.45\textwidth}
		\centering
		\includegraphics[width=\textwidth]{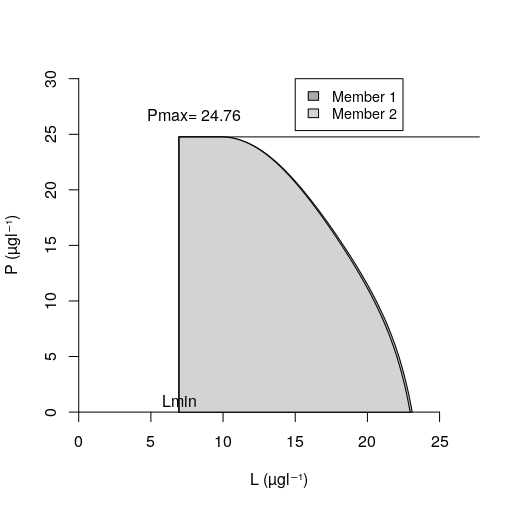}
		\caption{(a) Member 2 versus Member 1.}
		\label{fig:Member2and1}
	\end{subfigure}
	\hfill
	\begin{subfigure}[b]{0.45\textwidth}
		\centering
		\includegraphics[width=\textwidth]{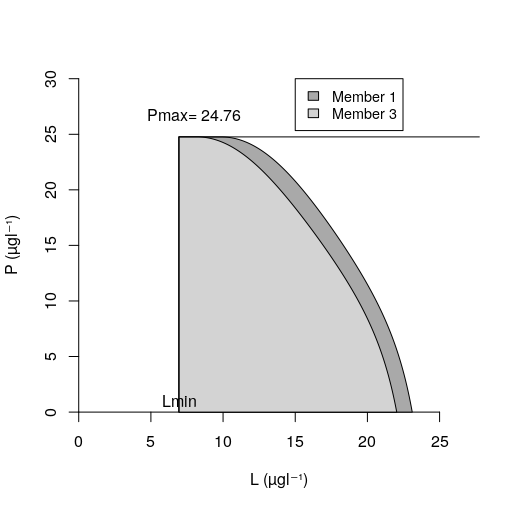}
		\caption{Member 3 versus Member 1.}
		\label{fig:Member3and1}
	\end{subfigure}
	\hfill
	\begin{subfigure}[b]{0.45\textwidth}
		\centering
		\includegraphics[width=\textwidth]{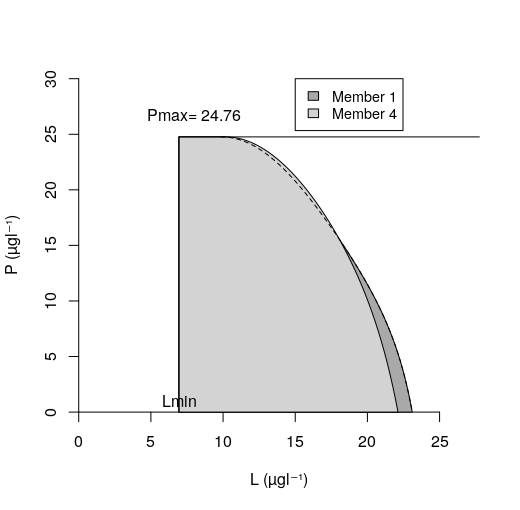}
		\caption{Member 4 guarantied viability kernel versus Member 1.}
		\label{fig:Member4and1}
	\end{subfigure}
	\hfill
	\begin{subfigure}[b]{0.45\textwidth}
		\centering
		\includegraphics[width=\textwidth]{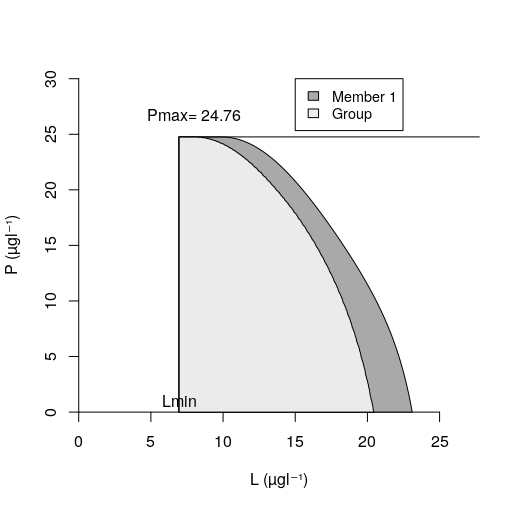}
		\caption{Guarantied Viability kernel $Guar_{f_{Bd}}$ associated to the group versus Member 1.}
		\label{fig:groupKernel}
	\end{subfigure}
	\caption{Guarantied viability kernel for each member and the group versus Member 1's viability kernel. Computation with R (\cite{R2010}) and ViabLab (\cite{viablab})}.
	\label{fig:Kernels}
\end{figure}

The committee members' believes we consider are summarized in Table \ref{tab:var}.
\begin{table}[b]
	\begin{tabular}{lcc@{\hskip.8mm}c@{\hskip.8mm}|c@{\hskip.8mm}c@{\hskip.8mm}c@{\hskip.8mm}c@{\hskip.8mm}}
		\hline
		Parameters & $b_i$ & $r_i$ & $m_i$	&  $\alpha_i$ & $q_i$ & $\lambda_i$ \\
	b, r, m 	  & P loss & max. & P value &  model & steepness &  steepness \\
	in $\mu g.l^{-1}$	 &  &rate & at half &type & model (\ref{eq:lake}) & model (\ref{eq:Bourget2}) \\
		 &  & &max. rate & & $S_i$  &  $S'_i$  \\

		\hline
Table \ref{tab:lake} &  \multirow{2}{*}{2.2676} &  \multirow{2}{*}{101.96} &  \multirow{2}{*}{26.90} & \multirow{2}{*}{1}  &  \multirow{2}{*}{2.222} & \multirow{2}{*}{-} \\
/ Member 1 &   &   &   &   &   &  \\
Member 2 & 2.2676 &  101.96 &  26.90 & 1  & $[2.2,2.3]$ & - \\
Member 3 & $[2.2,2.3]$ &  101.96 &  26.90 & 1  & 2.222 & - \\
Member 4 & 2.2676 &  101.96 &  26.90 & 0  & - & $[1/19,1/16]$ \\
		\hline
	\end{tabular}
	\caption{Believes regarding the model and parameters of Lake Bourget dynamics.}
	\label{tab:var}
\end{table}
Regarding the definition of the constraint set, we consider that agricultural activity leads to at least 25 tons of incoming P each year. This value is arbitrary but it is lower than the mean loading between 2004 and 2016 which was above 33 tons/year (see \cite{brias2018}). So we choose as lower limit $L_{min}=25/v\approx6.94 \; \mu g.l^{-1}$. Considering the desirable threshold for in-lake P, we consider an optimistic scenario with the value of the mesotropic equilibrium as maximum, $P_{max}=24.76 \; \mu g.l^{-1}$. 
The constraint set for this scenario is  $K=\left\lbrace (L,P), L\geq L_{min}, P\leq P_1\right\rbrace $ (respectively $ P_2$ for $K_2$). As possibility of control, we consider 
that the maximum rate for the reduction of incoming P is half the maximum difference $\Delta$ of loading between two consecutive years between 2004 and 2016. For the increase of the loading we consider that the maximum rate can be $\Delta$. So the set of admissible control is $U=[-\frac{\Delta}{2}, \Delta]$, with $\Delta\approx 3.15\; \mu g.l^{-1}y^{-1}$. 
When parameters are in a range, we consider the embedding dynamics $Sv$ as in equation \ref{ex1G} with the corresponding parameter as $v$ and $V$ its range.  
For each member $i$ 
it is possible to define a viability problem either as in Section \ref{sec:individual}, when parameters have fixed values, or a guarantied viability problem as in Section \ref{sec:GV}, when parameters are in a range. We note $\viab_i$ the viability kernel associated to Member i's project $(K,U)$ and dynamics $S_i$ (or $S'_i$) with fixed parameters (or $Sv_i$ or $S'_i$ with parameter in a range). Figure \ref{fig:Kernels} shows the viability and guarantied viability kernels computed for scenario 1 for members 2 to 4.

Applying the method described in the previous section, we define an embedding function $f_B$ for the group from model ()\ref{eq:Bourget2}) and Table \ref{tab:lake}. We then define the guarantied viability problem $Bv$ (equation \ref{eq:lakeBourgetV}) associated to $f_B$.
\begin{equation}
f_B(x,u,v)=\left( 
\begin{array}{l}
u\\
- v_1 P + L +  r \left( (1-v_2) \frac{P^{v_3}}{m^{v_3} + P^{v_3}} + v_2 \frac{P}{P + m e^{(-v_4(P-m)) }}\right)  \;
\end{array}
\right)
\label{eq:lakeBourget3}
\end{equation}
where $v_1$ stands for parameter $b_i$, $v_2$ for $\alpha_i$, $v_3$ for $q_i$ and $v_4$ for $\lambda_i$, with :
\begin{equation}
\left\{
\begin{array}{lcl}
x&=&(L,P) \in K\\
u&\in& U=[-\frac{\Delta}{2}, \Delta]\\
v &\in& V=[2.2,2.3]\times[0,1]\times[2.2,2.3]\times[1/19,1/18]\; 
\end{array}
\right.
\label{eq:paramlakeBourgetV}
\end{equation}
\begin{equation}
Bv\left\{
\begin{array}{lcl}
(L,P)'(t)&=&f_B((L,P)(t),u(t),v(t))\\
u(t)&\in &U\\
v(t)&\in &V \\
(L,P)(t) &\in& K \; 
\end{array}
\right.
\label{eq:lakeBourgetV}
\end{equation}
Since $U$ and $V$ are constant function of $(L,P)$, and since $f_B$ is Lipschitz, the conditions of Proposition \ref{prop:viabguaranty} are fulfilled. Since $U$ is constant, it is Lipschitz for every $\lambda>0$, so the guarantied viability kernel $Guar_\lambda,_{f_B}(K)$ associated to problem (\ref{eq:lakeBourgetV})  is closed and has the property of guarantied viability. It is then from Theorem \ref{th:viab_consensus} a consensus set of states for the four member of the committee whose believes are summarized in Table \ref{tab:var}. To compute an approximation of  $Guar_\lambda,_{f_B}(K)$, we use the ViabLab library \cite{viablab}, developed by  A. Désilles and used in \cite{durand2017agroecological}. This library uses the convergence conditions of the algorithm established by P. Saint-Pierre \cite{saint1994approximation}. Since the ViabLab library requires presently for the computation of guarantied viability kernel discrete problems in time and space, we defined a discretized version of the viability problem, with function $f_{Bd}$ from equation (\ref{eq:lakeBourget4}), with a discretization parameter $\tau=0.1$ for which the dynamics are stable. We also used projection on grid method from \cite{Lavallee2020v2} to minimize discretization error.  
\begin{equation}
f_{Bd}((L,P),u,v)=\left( 
\begin{array}{l}
L + \tau u\\
P + \tau \left[ - v_1 P + L + \right.   \\
\left. \qquad r \left(  (1-v_2) \frac{P^{v_3}}{m^{v_3} + P^{v_3}} + v_2 \frac{P}{x_2 + m e^{(-v_4(P-m)) }}\right) \right]   \;
\end{array}
\right)
\label{eq:lakeBourget4}
\end{equation}
The resulting approximation $Guar_{f_{Bd}}(K)$ is shown on Figure \ref{fig:groupKernel}. The guarantied viability kernel for the group is a viable set for all members, and the viable controls on its boundary are viable controls for all members.

%
%
\subsection{Discussion}
\label{sec:discuss}
Depending on the dynamics and the beliefs of the different group members, the guarantied viability kernel computed following the approach of Section \ref{sec:PEF} could be smaller than the largest viable set corresponding to the union of the parameters set of each group member. For instance, with $(b,r) \in \{(2.1,100),(2.2,80)\}$, it is possible to design a viable set for these two values only, solving the problem with the computation of integral curves as it is done in Section \ref{sec:viewpoints} for the lake problem (see figure \ref{fig:celacgaranti}). Whereas following the method in Section \ref{sec:PEF}, in order to respect the conditions of VT theorems and use the viablab library it is necessary to define a guaranty viability problem for $(b,r) \in [2.1,2.2]\times[80,100]$. But when the dimension of the state space is greater than 2, the first method is virtually impossible to implement with a generic module (since a specific mathematical study of the dynamics is necessary). 

The viability algorithm is exponential with the dimension of the space in the general case, so it can be very slow, in particular when the viability kernel is empty or with high dimension problems. For each alternative model definition a tyche has to be considered, which increases the dimension of $V$ linearly with the number of models.
The computation of $Guar_{f_{Bd}}(K)$ on Figure \ref{fig:groupKernel} takes 938.64s with a processor i7-8650U CPU @ 1.90GHz × 8 and 15.5GiB RAM, with an accuracy of 1000 points/axis and a discretization step of 11 for control and the model type and 5 for the three other tyches ($b$,$q$,$\lambda$). It takes $1243.86\; s$ for a scenario with $P_{max} = 15.0  \; \mu g.l^{-1}$  (all parameters being equal), and the guarantied viability kernel is not empty. Increasing the dimension of $V$ can also lead to much longer computation time. For instance, it takes $7122.99\; s$ for the scenario with $P_{max} = 15.0  \; \mu g.l^{-1}$, all parameters being equal except for $m \in [26.0,27.0]$ as additional tyche. In that case the guarantied viability kernel is empty. This method has been used for a dimension 3 problem of management of marine protected area \cite{Zaleski2020}. With an accuracy of 100 points/axis and dimension 2 controls and tyches (with discretization step of 11), the computation time is $437.94 \; s$. It is $3153.31 \; s$ with an accuracy of 300 points/axis and 5 steps for each tyche.

Regarding the model of Lake Bourget itself, we consider here a single control for different practices (use of wetlands, use of retention basins, different farmer practices), and the value of its range is  consistent with observations but arbitrary. The model could be improved by taking into account more detailed mechanisms for these different types of control and their relation to soil leaching and rainfall. 
The lake is considered as homogeneous, which seams a reasonable assumption regarding the water resident time of 14 years (\cite{brias2018}). On the other hand, since blooms of cyanobacteria are often localized, it could be useful to use spatial and weekly data to assess the size of perturbations and take them into account to consider robustness issues as defined in \cite{robust2019}.  




\section{Conclusion}
In this paper, we have presented a method for reaching a viability-based consensus for the management of commons with multiple usages. We have proposed a definition of the management project and stakeholder viewpoints, and a viability-based definition for the consensus solution as a viable set for all stakeholders with conditions on their regulation map. We have defined embedding functions that allow to compute a guarantied viability kernel for the associated dynamics. We have then shown that this guarantied viability kernel is a consensus solution (Theorem \ref{th:viab_consensus}). It can then be computed with algorithms used for viability kernel approximation. We have then applied this method to a management scenario for Lake Bourget. The main interest of this method is that stakeholders can retain their vision of the dynamics. Negotiations can focus on the definition of desirable states and admissible actions. It also prepare the way for an alternative to agent based modeling when dealing with stakeholders for management of commons.


\bibliographystyle{apalike}        
\bibliography{collectiveGestionLakeEM}           

\begin{thebibliography}{}

\bibitem[Alais et~al., 2017]{Alais_DeLara_multi-usage_2017}
Alais, J.-C., Carpentier, P., and De~Lara, M. (2017).
\newblock Multi-usage hydropower single dam management: chance-constrained
  optimization and stochastic viability.
\newblock {\em Energy Systems}, 8(1):7--30.

\bibitem[Aubin, 1991]{Aubin91}
Aubin, J.-P. (1991).
\newblock {\em Viability Theory}.
\newblock {Birkh\"{a}user}, Basel.

\bibitem[Aubin, 1997]{AubinGuaranti}
Aubin, J.-P. (1997).
\newblock {\em Dynamic economic theory: a viability approach}, volume~5.
\newblock Springer Verlag.

\bibitem[Barreteau et~al., 2001]{o_barreteau_role-playing_nodate}
Barreteau, O., Bousquet, F., and Attonaty, J.-M. (2001).
\newblock Role-playing games for opening the black box of multi-agent systems:
  method and lessons of its application to {Senegal} {River}.
\newblock {\em Journal of Artificial Societies and Social Simulation}, 4(2):5.

\bibitem[Bates and Saint-Pierre, 2018]{BATES2018244}
Bates, S. and Saint-Pierre, P. (2018).
\newblock Adaptive policy framework through the lens of the viability theory: A
  theoretical contribution to sustainability in the anthropocene era.
\newblock {\em Ecological Economics}, 145:244--262.

\bibitem[Brias et~al., 2018]{brias2018}
Brias, A., Mathias, J.-D., and Deffuant, G. (2018).
\newblock Inter-annual rainfall variability may foster lake regime shifts: {An}
  example from {Lake} {Bourget} in {France}.
\newblock {\em Ecological Modelling}, 389:11--18.

\bibitem[Carpenter et~al., 1999]{carpenter_management_1999}
Carpenter, S.~R., Ludwig, D., and Brock, W.~A. (1999).
\newblock Management of {Eutrophication} for {Lakes} {Subject} to {Potentially}
  {Irreversible} {Change}.
\newblock {\em Ecological Applications}, 9(3):751--771.

\bibitem[{de Groot}, 2006]{deGroot2005}
{de Groot}, R. (2006).
\newblock Function-analysis and valuation as a tool to assess land use
  conflicts in planning for sustainable, multi-functional landscapes.
\newblock {\em Landscape and Urban Planning}, 75(3):175 -- 186.

\bibitem[De~Lara and Martinet, 2009]{23_de_lara_multi-criteria_2009}
De~Lara, M. and Martinet, V. (2009).
\newblock Multi-criteria dynamic decision under uncertainty: {A} stochastic
  viability analysis and an application to sustainable fishery management.
\newblock {\em Mathematical Biosciences}, 217(2):118--124.

\bibitem[D{\'e}silles, 2020]{viablab}
D{\'e}silles, A. (2020).
\newblock Viablab library.
\newblock Technical report, LASTRES. https://github.com/lastre-viab/VIABLAB.

\bibitem[Dowling et~al., 2020]{DOWLING2020109243}
Dowling, N.~A., Dichmont, C.~M., Leigh, G.~M., Pascoe, S., Pears, R.~J.,
  Roberts, T., Breen, S., Cannard, T., Mamula, A., and Mangel, M. (2020).
\newblock Optimising harvest strategies over multiple objectives and
  stakeholder preferences.
\newblock {\em Ecological Modelling}, 435:109243.

\bibitem[Doyen, 2000]{doyen_lipschitz_2000}
Doyen, L. (2000).
\newblock Lipschitz {Kernel} of a {Closed} {Set}-{Valued} {Map}.
\newblock {\em Set-Valued Analysis}, 8(1):101--109.

\bibitem[Durand et~al., 2017]{durand2017agroecological}
Durand, M.-H., D{\'e}silles, A., Saint-Pierre, P., Angeon, V., and
  Ozier-Lafontaine, H. (2017).
\newblock Agroecological transition: A viability model to assess soil
  restoration.
\newblock {\em Natural resource modeling}, 30(3):e12134.

\bibitem[Etienne, 2014]{ComMod_2014}
Etienne, M., editor (2014).
\newblock {\em Companion {Modelling}: {A} {Participatory} {Approach} to
  {Support} {Sustainable} {Development}}.
\newblock Springer Netherlands.

\bibitem[Gajardo et~al., 2017]{gajardo_modeling_2017}
Gajardo, P., Harmand, J., Ramirez, H., Rapaport, A., Riquelme, V., and
  Rousseau, A. (2017).
\newblock Modeling and control of in-situ decontamination of large water
  resources.
\newblock {\em ESAIM: Proceedings and Surveys}, 57:70--85.

\bibitem[Hardin, 1968]{hardin1968}
Hardin, G. (1968).
\newblock The {Tragedy} of the {Commons}.
\newblock {\em Science, New Series}, 162(3859):1243--1248.

\bibitem[Jacquet, 2018]{Jacquet2017}
Jacquet, S. (2018).
\newblock {Suivi scientifique du lac du Bourget pour l'ann{\'e}e 2017}.
\newblock Technical report, {INRAE ; CISALB ; Agence de l'eau RMC}.

\bibitem[Křivan, 1991]{53_krivan_differential_1991}
Křivan, V. (1991).
\newblock Differential inclusions as a methodology tool in population biology.
\newblock In {\em Proceedings of the 1995 {European} {Simulation}
  {Multiconference}}, pages 544--547, Prague.

\bibitem[Křivan and Colombo, 1998]{55_krivan_non-stochastic_1998}
Křivan, V. and Colombo, G. (1998).
\newblock A {Non}-stochastic {Approach} for {Modeling} {Uncertainty} in
  {Population} {Dynamics}.
\newblock {\em Bulletin of Mathematical Biology}, 60(4):721--751.

\bibitem[Lavall{\'e}e, 2020]{Lavallee2020v2}
Lavall{\'e}e, F. (2020).
\newblock PhD thesis.
\newblock Model, exploration, viability theory for management support of
  invasive species. Sorbonne University. In French.

\bibitem[Little et~al., 2007]{44_1_little_elfsimmodel_2007}
Little, L.~R., Punt, A.~E., Mapstone, B.~D., Pantus, F., Smith, A. D.~M.,
  Davies, C.~R., and McDonald, A.~D. (2007).
\newblock {ELFSim}—{A} model for evaluating management options for spatially
  structured reef fish populations.
\newblock {\em Ecological Modelling}, 205(3):381--396.

\bibitem[Liu et~al., 2015]{LakeChinois}
Liu, W., Wang, S., Zhang, L., and Ni, Z. (2015).
\newblock Water pollution characteristics of {Dianchi} {Lake} and the course of
  protection and pollution management.
\newblock {\em Environmental Earth Sciences}, 74(5):3767--3780.

\bibitem[Martin, 2004]{Martin2004}
Martin, S. (2004).
\newblock The cost of restoration as a way of defining resilience: a viability
  approach applied to a model of lake eutrophication.
\newblock {\em Ecol. Soc.}, 2(9):8.

\bibitem[Martin and Alvarez, 2019]{robust2019}
Martin, S. and Alvarez, I. (2019).
\newblock Anticipating shocks in the state space: Characterizing robustness and
  building increasingly robust evolutions.
\newblock {\em SIAM Journal on Control and Optimization}, 57(1):490--509.

\bibitem[Martinet et~al., 2016]{45_martinet_risk_2016}
Martinet, V., Peña-Torres, J., De~Lara, M., and Ramírez~C., H. (2016).
\newblock Risk and {Sustainability}: {Assessing} {Fishery} {Management}
  {Strategies}.
\newblock {\em Environmental and Resource Economics}, 64(4):683--707.

\bibitem[Oubraham and Zaccour, 2018]{zaccour_survey_2018}
Oubraham, A. and Zaccour, G. (2018).
\newblock A {Survey} of {Applications} of {Viability} {Theory} to the
  {Sustainable} {Exploitation} of {Renewable} {Resources}.
\newblock {\em Ecological Economics}, 145:346--367.

\bibitem[Parrachino et~al., 2006]{parrachino_cooperative_2006}
Parrachino, I., Dinar, A., and Patrone, F. (2006).
\newblock {\em Cooperative {Game} {Theory} {And} {Its} {Application} {To}
  {Natural}, {Environmental}, {And} {Water} {Resource} {Issues} : 3.
  {Application} {To} {Water} {Resources}}.
\newblock Policy {Research} {Working} {Papers}. The World Bank.

\bibitem[Saint-Pierre, 1994]{saint1994approximation}
Saint-Pierre, P. (1994).
\newblock Approximation of the viability kernel.
\newblock {\em Applied Mathematics and Optimization}, 29(2):187--209.

\bibitem[Sanogo et~al., 2012]{38_sanogo_viability_2012}
Sanogo, C., Ben~Miled, S., and Raissi, N. (2012).
\newblock Viability {Analysis} of {Multi}-fishery.
\newblock {\em Acta Biotheoretica}, 60(1):189--207.

\bibitem[Team, 2010]{R2010}
Team, R. D.~C. (2010).
\newblock {\em R: A language and environment for statistical computing}.
\newblock {R Foundation for Statistical Computing}, Vienna.

\bibitem[Vinçon-Leite and Tassin, 1990]{vincon-leite_modelisation_1990}
Vinçon-Leite, B. and Tassin, B. (1990).
\newblock Modélisation de la qualité des eaux des lacs profonds : modèle
  thermique et biogéochimique du lac du {Bourget}.
\newblock {\em La Houille Blanche}, (3-4):231--236.

\bibitem[Vollenweider, 1982]{Vollenweider}
Vollenweider (1982).
\newblock coord. {E}utrophication of waters. {M}onitoring, assessment and
  control.
\newblock Technical report, OECD.

\bibitem[Wei et~al., 2013]{wei_sustainability_2013}
Wei, W., Alvarez, I., and Martin, S. (2013).
\newblock Sustainability analysis: {Viability} concepts to consider transient
  and asymptotical dynamics in socio-ecological tourism-based systems.
\newblock {\em Ecological Modelling}, 251:103--113.

\bibitem[Zaleski, 2020]{Zaleski2020}
Zaleski, L. (2020).
\newblock PhD thesis.
\newblock Decision and negociation assistant by viability analysis: Application
  to the participatory management of protected areas. Sorbonne University. In
  French.

\end{thebibliography}
\appendix
\section{Properties of Viability Kernels}
\label{annexe1}
\subsection{Properties of Multi-valued Maps}
Let $G$ be a multi-valued map from $\R^n$ to $\R^p$. The domain of $G$ is $Dom(G)=\left\lbrace x \in \R^n, G(x) \neq \emptyset\right\rbrace $. The graph of $G$ is $Graph(G)=\{ (x,y) \in \R^n \times \R^p, y \in G(x) \} $. $G$ has a linear growth if there is $c>0$ such that for all $x \in Dom(G)$, $\lVert G(x)\rVert \leq c(\lVert x \rVert +1)$. The system $Sc(f,U)$ (\ref{evolve}) is Marchaud if $f$ is continuous, $Graph(U)$ is closed, $f$ and $U$ have linear growth and the image set $\{f(x,u), u \in U(x)\} $ is convex for all of $x \in Dom(U)$. $G$ is Lipschitz for constant $\lambda>0$ (or $\lambda$-Lipschitz) if for all $x_1$, $x_2$ in $\R^n$, $G(x_1) \subset G(x_2) + \lambda \lVert x_1 - x_2 \rVert B$, where B is the unit ball.

\subsection{Closed Viability Kernels \cite{ Aubin91}}
\begin{prop}
	Continuous case: When the system $Sc(f,U)$ (\ref{evolve}) is Marchaud and $K$ is closed, the associated viability kernel is closed. It is the largest viable set in $K$.
	\label{prop:kernelclosed_c}
\end{prop}

\begin{prop}
	Discrete case: When system $Sd(f,U)$ (\ref{evolveD}) is such that $f$ is continuous, $U$ has a linear growth, $Graph(U)$ is closed, and when $K$ is closed, the associated viability kernel is closed.
	\label{prop:kernelclosed_d}
\end{prop}

\subsection{Guaranteed Viability Kernels}
\label{annexe1.2}
We recall some properties of dynamical controlled tychastic systems. 
\begin{mydef}
	A dynamical controlled tychastic system $Scv(f,U,V)$ (\ref{evolveV}) is Lipschitz if $f$ is Lipschitz, and $U$ and $V$ are Lipschitz with compact images. 
\end{mydef}
\begin{prop}
	From \cite{doyen_lipschitz_2000} in the continuous case and \cite{Lavallee2020v2} in the discrete case. The guaranteed viability kernel associated to a set $K$ is closed when the dynamics verify the following conditions: In the continuous case, when $K$ is closed and $Scv(f,U,V)$ is Lipschitz ; In the discrete case, when $f$ and $V$ are continuous, $Graph(U)$ is closed and $U$ has a linear growth (see \ref{sec:Viability} for the definitions). 
	\label{prop:viabguaranty}
\end{prop}
We recall Theorem \ref{th:viab_consensus}: 	The guaranteed viability kernel associated to $K$ for $Scv(f,U,V)$ (with $\lambda$ Lipschitz constant in the continuous case) (resp. $Sdv(f,U,V)$ in the discrete case ) is a consensus solution to the management of $\A$.

\textbf{Proof of Theorem \ref{th:viab_consensus}}.
Let $L\neq \emptyset$ be the guaranteed viability kernel for system (\ref{evolveV}) with $\lambda$ Lipschitz constant (resp. system (\ref{evolveDV}) in the discrete case) associated to constraint set $K$. Let $\tilde{U}$ be the guaranteed regulation map. By definition, $Dom\tilde{U}=L$. We now prove that the guaranteed viability kernel is a viable set for each member $i$. Let $i \in \N$, we consider the system $Sc'_i=Sc(f_i,\tilde{U})$ (resp. $Sd'_i=Sd(f_i,\tilde{U})$ in the discrete case). Since the guaranteed viability kernel is defined from the group project $(K,U)$, we have $L \subset K$, and for all $x\in L$, $\tilde{U}(x) \subset U(x)$. So an evolution governed by $Sc'_i$ (resp. $Sd'_i$) is also an evolution governed by $Sc(f_i,U)$ (resp. $Sd(f_i,U)$). 
Since system (\ref{evolveV}) (resp. system (\ref{evolveDV}) ) embeds $Sc(f_i,U_i)$ (resp. $Sd(f_i,U_i)$), it also embeds $Sc'_i$ (resp. $Sd'_i$). Let $x_0 \in L$, and let $x(.)$ (resp. $x^k$) be a trajectory starting at $x_0$ and governed by $Sc'_i$ (resp. $Sd'_i$). Because of the embedding there is a function $v$ such that $f_i(x(t),\tilde{u}(t)) = f(x(t),\tilde{u}(t),v(t))$ (resp. $f_i(x^k),\tilde{u}^k) = f(x^k,\tilde{u}^k,v^k)$ in the discrete case). So $x(.)$ is also an evolution governed by system (\ref{evolveV}) (resp.(\ref{evolveDV})). Since $L$ is the guaranteed viability kernel for system (\ref{evolveV}) with $\lambda$ Lipschitz constant (resp.(\ref{evolveDV})) associated to constraint $K$, from Definition (\ref{guarantiSet}), all trajectories starting from $x_0 \in L$ and governed by (\ref{evolveV}) (resp.(\ref{evolveDV})) are viable in $L$ for control selection in $\tilde{U}$. So the trajectory of $x(.)$ governed by $Sc'_i$ (resp. $Sd'_i$) starting at $x_0$ stays in $L$. So $x(.)$ is an evolution starting at $x_0$  governed by $Sc_i$ (resp. $Sd_i$) viable in $L$. So $L$ is a viable set for member $i$. 
$\blacksquare$
\section{Embedding system}
\label{annexe2}
\begin{prop}
	If $U$ has a linear growth and for all $i \in \N$, $f_i$ has a linear growth, then a System (\ref{evolveV}) (resp. \ref{evolveDV}) can embed $Sc(f_i,U)$(\ref{evolve}) (resp. $Sd(f_i,U)$ (\ref{evolveD})) for all $i \in \N$.
\end{prop} 
\textbf{Proof.} We note $M_x = \max_{i\in \N, u\in U(x)}(\lVert f_1(x,u)-f_i(x,u)\rVert)$. $M_x$ is defined since $U$ and all $f_i$ have a linear growth. We define $f(x,u,v)=f_1(x,u)+v(x)$ with $v(x) \in  V(x)=B(0,M_x)$, where $B(a,r)$ is the closed ball with center $a$ and radius $r$.
Then  for $x \in K$ and $u \in U(x)$ we define $v_{i,u,x} = f_i(x,u)-f_1(x,u)$. We have $\lVert v_{i,u,x}\rVert \leq M_x$, so $v_{i,u,x} \in V(x)$.
Then $\forall i \in \N, f_i(x,u) = f(x,u,v_{i,u,x}) $ and equation (\ref{condV}) is verified.
$\blacksquare$

Definition of $f$ leading to smaller sets of perturbation are preferable. For instance, it can be interesting to define $f$ with the convex hull of the $f_i(x,u)$: With $i \in J=\N \smallsetminus \{1\}$ we consider $v=(v_i), v_i \in [0,1]$ with $\sum_{i \in J} v_i\leq1$,  and define $f(x,u,v)= f_1(x,u)(1-\sum\limits_{i\in J} v_i) + \sum\limits_{i\in J}   v_i f_i(x,u)$.
\end{document}